\documentclass[12pt]{svmult}
\usepackage{mathptmx, helvet, courier, graphicx, makeidx, multicol, footmisc, amsmath,amsfonts,amssymb}

\newtheorem{defn}[theorem]{Definition}

\def\irr#1{{\rm Irr}(#1)}
\def\irrr#1#2 {\irr {#1 \mid #2}}
\newcommand{\R}{\mathbb R}
\newcommand{\n}{\emph{n}}
\newcommand{\C}{\mathbb [(n+1)\nu_{n+1}J_n]^{-1}}
\newcommand{\F}{e^{-\varphi(|y|)}}

\begin{document}

\title*{Maximal surface area of a convex set in $\R^n$ with respect to log concave rotation invariant measures.}
\author{Galyna Livshyts}

\institute{Galyna Livshyts \at Department of Mathematical Sciences,  Kent State University,  Kent, OH, 44242, USA \email{glivshyt@kent.edu}}
\titlerunning{Surface area with respect to log concave rotation invariant measures.}
\authorrunning{Galyna Livshyts}
\maketitle

\abstract*{It was shown by K. Ball and F. Nazarov, that the maximal surface area of a convex set in $\R^n$ with respect to the Standard Gaussian measure is of order
$n^{\frac{1}{4}}$. In the present paper we establish the analogous result for all rotation invariant log concave probability measures. We show that the maximal surface area
with respect to such measures is of order $\frac{\sqrt{n}}{\sqrt[4]{Var|X|} \sqrt{\mathbb{E}|X|}}$, where $X$ is a random vector in $\R^n$ distributed with respect to the measure.
}
\section {Introduction}

In this paper we will study geometric properties of the probability measures  $\gamma$ on $\R^n$ with density $C_{n} e^{-\varphi(|y|)},$ where $\varphi(t)$ is a nonnegative nondecreasing convex function, which may take infinity as a value, and the normalizing constant $$C_n=\left(\int_{\R^n}e^{-\varphi(|y|)}dy\right)^{-1}.$$

We recall that the Minkowski surface area of a convex set $Q$ with respect to the measure $\gamma$ is defined to be
\begin{equation}\label{outer}
\gamma(\partial Q)=\liminf_{\varepsilon\rightarrow +0}\frac{\gamma((Q+\varepsilon B_2^n)\backslash Q)}{\varepsilon},
\end{equation}
where $B_2^n$ denotes the Euclidian unit ball in $\R^n.$

The special case of $\varphi(t)=\frac{t^2}{2}$, which corresponds to the standard Gaussian measure $\gamma_2$, has been actively studied. Sudakov, Tsirelson \cite{ST} and Borell \cite{B} proved, that among all convex sets of a fixed Gaussian measure, half spaces have the smallest Gaussian surface area. Mushtari and Kwapien asked the reverse version of isoperimetric inequality, i.e. how large the Gaussian surface area of a convex set $A\subset \R^n$ can be. It was shown by Ball \cite{ball}, that Gaussian surface area of a convex set in $\R^n$ is asymptotically bounded by $C n^{\frac{1}{4}}$, where $C$ is an absolute constant. Nazarov \cite{fedia} proved the sharpness of Ball's result and gave the complete solution to this asymptotic problem:
\begin{equation}\label{fedia}
0.28 n^{\frac{1}{4}}\leq \max_{Q\in \mathcal{K}_n} \gamma_2(\partial Q)\leq 0.64 n^{\frac{1}{4}},
\end{equation}
where by $\mathcal{K}_n$ we denote the set of all convex sets in $\R^n$. Further estimates for $\gamma_2(\partial Q)$ for the special case of polynomial level set surfaces were provided by D. Kane \cite{kane}. He showed that for any polynomial $P(y)$ of degree $d$, $\gamma_2(P(y)=0)\leq \frac{d}{\sqrt{2}}$.

Isoperimetric inequalities for a wider class of rotation invariant measures were studied by Sudakov and Tsirelson \cite{ST}. Recently, geometric properties for various classes of rotation invariant measures were established by Bobkov \cite{bob1, bob2, bob3}, Bray and Morgan \cite{BM}, Maurmann and Morgan \cite{MM} and others.

The maximal surface area of convex sets for the probability measures $\gamma_p$ with densities $C_{n,p}e^{-\frac{|x|^p}{p}}$, where $p>0$, was studied in \cite{L}. It was shown there, that
\begin{equation}\label{ya}
c(p)n^{\frac{3}{4}-\frac{1}{p}}\leq \max_{Q\in \mathcal{K}_n} \gamma_p(\partial Q) \le C(p)n^{\frac{3}{4}-\frac{1}{p}},
\end{equation}
where $c(p)$ and $C(p)$ are constants depending on $p$ only.

In the present paper we obtain a generalization of results due to Ball and Nazarov, and find an expression for the maximal surface area with respect to an arbitrary rotation invariant log concave measure $\gamma$. The expression depends on the measure's natural characteristics, i.e. expectation and variance of a random variable, distributed with respect to $\gamma.$

We shall use notation $\precsim$ for an asymptotic inequality: we say that $A(n)\precsim B(n)$ if there exists an absolute positive constant $C$ (independent of $n$), such that $A(n)\leq CB(n)$. Correspondingly, $A(n)\approx B(n)$ means that $B(n)\precsim A(n)\precsim B(n)$.


The following theorem is the main result of the present paper:
\begin{theorem}\label{th:main}
Fix $n\geq 2$. Let $\gamma$ be log concave rotation invariant measure on $\R^n.$ Consider a random vector $X$ in $\R^{n}$ distributed with respect to $\gamma$. Then
$$\max_{Q\in \mathcal{K}_n}\gamma(\partial Q)\approx\frac{\sqrt{n}}{\sqrt{\mathbb{E}|X|} \sqrt[4]{Var|X|}},$$
where, as usual, $\mathbb{E}|X|$ and $Var|X|$ denote the expectation and the variance of $|X|$ correspondingly.
\end{theorem}

Let us note, that the above Theorem implies (\ref{fedia}). It also implies (\ref{ya}) in the case $p\geq 1,$
and the details of these implications are shown in Section 3.

Another classical example of a log concave rotation invariant measure is the normalized Lebesgue
measure restricted to the unit ball. In that case $\varphi(t)$ equals to zero for all $t<1$ and takes infinity as a value for all $t\geq 1$.
For that measure $\mathbb{E}|X|\approx 1$ and $Var|X|\approx n^{-2}$, so the maximal surface area is of order
$\frac{\sqrt{n}}{\sqrt[4]{n^{-2}\times 1}}=n$. The set with the maximal surface area is the sphere of radius $1$, which is also clear by monotonicity of the standard surface area measure. 

We outline, that for isotropic measures (see \cite{MP} for definitions and details), Theorem \ref{th:main} together with the result from \cite{Kl2} entails that the maximal perimeter varies between $C_1 n^{\frac{1}{4}}$ and $C_2 n^{\frac{1}{2}}$, where $C_1$ and $C_2$ are absolute constants. The standard Gaussian measure is an example of an isotropic measure with the maximal surface area of order $n^{\frac{1}{4}}$, and the Lebesgue measure restricted to a ball of radius $\sqrt{n}$ is an example of an isotropic measure with the maximal surface area of order $n^{\frac{1}{2}}$. If the measures $\gamma_p$ from (\ref{ya}) are brought to the isotropic position, the maximal surface area with respect to them is of order $n^{\frac{1}{4}}$.

The main definitions, technical lemmas and some preliminary facts are given in Section 2.
Some connections between the probabilistic and analytic setup are provided in Section 3. The upper bound for  Theorem \ref{th:main} is obtained in
Section 4, and the lower bound is shown in Section 5. In Section 6 we provide some examples to exhibit the sharpness
of Theorem \ref{th:main}.

\noindent {\bf Acknowledgment}. I would like to thank Artem Zvavitch and Fedor Nazarov for introducing me to the subject, suggesting me this problem and for extremely helpful and fruitful discussions. I would also like to thank Benjamin Jaye for a number of useful remarks. The author is supported in part by U.S. National Science Foundation grant DMS-1101636.

\section{Some definitions and lemmas.}

 This section is dedicated to some general properties of spherically invariant log concave measures.
 We outline some basic facts which are needed for the proof. Some of them have appeared in literature.
 See \cite{Kl} for an excellent overview of the properties of log concave measures, in particular the proof of Lemma 4.5, where some portion of the current section
 appears.

 We write all the calculations in $\R^{n+1}$ instead of $\R^n$ for the notational simplicity. We use notation $|\cdot|$ for the norm in Euclidean space $\R^{n+1}$;
  $|A|$ stands for the Lebesgue measure  of a measurable set $A \subset \R^{n+1}$.
 We will write     $B_2^{n+1}=\{x\in \R^{n+1}: |x| \le 1\}$ for the unit ball in $\R^{n+1}$ and $\mathbb{S}^{n}=\{x\in \R^{n+1}: |x|=1\}$ for the
 unit sphere. We denote $\nu_{n+1}=|B_2^{n+1}|=\frac{\pi^{\frac{n+1}{2}}}{\Gamma(\frac{n+1}{2}+1)}$. We note that $|\mathbb{S}^n|=(n+1)\nu_{n+1}$.

We fix a convex nondecreasing function $\varphi(t):[0,\infty)\rightarrow [0,\infty]$. Let $\gamma$ be a probability measure on $\R^{n+1}$ with density $C_{n+1}\F$. The normalizing constant $C_{n+1}$ equals to $[(n+1)\nu_{n+1} J_{n}]^{-1}$, where
\begin{equation}\label{maint}
J_{n}=\int_0^{\infty} t^{n} e^{-\varphi(t)} dt.
\end{equation}
We introduce the notation $g_n(t)=t^n e^{-\varphi(t)}$. Since we normalize the measure anyway, we may assume that $\varphi(0)=0$.

Without loss of generality we may assume that $\varphi\in C^2[0,\infty)$. This can be shown by the standard smoothing argument (see, for example, \cite{cord}).

We shall use a well known integral formula for $\gamma(\partial Q)$, which holds true, in particular, for the measures with continuous densities:
\begin{equation}\label{def2}
\gamma(\partial Q)=C_{n+1}\int_{\partial Q} \F d\sigma(y),
\end{equation}
where $d\sigma(y)$ stands for Lebesgue surface measure (see the Appendix for the proof). 

The below Lemma shows that the surface area with respect to $\gamma$ is stable under small perturbations.

\begin{lemma}\label{homothety}
Fix $n\geq 2$. Let $M$ be a measurable subset of a boundary of a convex set in $\R^{n+1}$. Then
$$C_{n+1}\int_{\frac{1}{1+\frac{1}{n}}M} e^{-\varphi(|y|)} d\sigma(y)\succsim \gamma(M).$$
\end{lemma}
\noindent\textbf{Proof.} We observe, that $\varphi\left(\frac{|y|}{1+\frac{1}{n}}\right)\leq \varphi(|y|)$, since $\varphi(t)$ is nondecreasing. Also, $$d\sigma\left(\frac{y}{1+\frac{1}{n}}\right)=\left(1+\frac{1}{n}\right)^{-n}d\sigma(y)\approx d\sigma(y).$$
We conclude:
$$
\int_{\frac{1}{1+\frac{1}{n}}M} e^{-\varphi(|y|)} d\sigma(y)=\int_{M} e^{-\varphi\left(\frac{|y|}{1+\frac{1}{n}}\right)} d\sigma\left(\frac{y}{1+\frac{1}{n}}\right)\succsim \gamma(M). \,\,\square
$$
\begin{remark}
We observe as well, that the same statement holds for all measures with densities, decreasing along each ray starting at zero. 
\end{remark}



\begin{defn}\label{tzerodef}
We define $t_0$ to be the point of maxima of the function $g_n(t)$, i.e., $t_0$ is the solution of the equation
\begin{equation}\label{todef}
\varphi'(t)t=n.
\end{equation}
\end{defn}
We note that the equation (\ref{todef}) has a solution, since $t\varphi'(t)$ is nondecreasing, continuous and $\lim_{t\rightarrow +\infty} t\varphi'(t)=+\infty$.
This solution is unique, since $t\varphi'(t)$ strictly increases on its support. This definition appears in most of the literature dedicated to
rotation invariant log concave measures: see, for example, \cite{KM} (Section 2), \cite{bob2} (Remark 3.4) or \cite{Kl} (Lemma 4.3).
\begin{remark}\label{to}
We may define $t_n$ and $t_{n-1}$ by
\begin{equation}\label{tnmin1}
\varphi'(t_{n-1})t_{n-1}=n-1,
\end{equation}
\begin{equation}\label{tn}
\varphi'(t_n)t_n=n.
\end{equation}
We claim that $t_n^n\approx t_{n-1}^n$. To see this, we note that function $t\varphi'(t)$ is nondecreasing.
Hence $t_n\geq t_{n-1}$. On the other hand, subtracting (\ref{tnmin1}) from (\ref{tn}), we get
$$1=\varphi'(t_n)t_n-\varphi'(t_{n-1})t_{n-1}\geq \varphi'(t_{n-1}) (t_n-t_{n-1})=\frac{n-1}{t_{n-1}} (t_n-t_{n-1}).$$
The above leads to the following chain of inequalities:
\begin{equation}\label{compar}
1\leq\frac{t_{n}}{t_{n-1}}\leq 1+\frac{1}{n-1},
\end{equation}
and therefore $t_n^n\approx t_{n-1}^n$.

In a view of the above we introduce the notation ``$t_0$''. Everywhere in the paper it is assumed that $t_0=t_{n+1}$.

We notice in addition, that
\begin{equation}\label{philessthann}
\varphi(t_0)=\varphi(t_0)-\varphi(0)\leq \varphi'(t_0)t_0=n,
\end{equation}
since $\varphi(0)=0$ by our assumption.
\end{remark}

The next lemma provides simple asymptotic bounds for $J_n$. It was proved in \cite{KM}, but for the sake of completeness we sketch the proof below.
\begin{lemma}\label{mainintegral}
$$\frac{g_n(t_0)t_0}{n+1} \leq J_n\leq \sqrt{2\pi}(1+o(1))\frac{g_n(t_0)t_0}{\sqrt{n}}.$$
\end{lemma}
\noindent\textbf{Sketch of the proof.} The integral $J_n$ can be estimated from above by Laplace method, which can be found, for example, in \cite{bruign}. We rewrite
\begin{equation}\label{jn}
J_n=g_n(t_0)\int_0^{\infty} e^{n\log\frac{t}{t_0}+\varphi(t_0)-\varphi(t)}dt.
\end{equation}
By the Mean Value theorem, $\varphi(t_0)-\varphi(t)\leq \varphi'(t_0)(t_0-t)=n(1-\frac{t}{t_0})$ for any $t\geq 0$. Thus, (\ref{jn}) is less than
$$g_n(t_0)t_0\int_0^{\infty} e^{n h(t)}dt,$$
where $h(t)=\log t-t+1$.  It is easy to check that $h(t)$ satisfies Laplace's condition (see \cite{bruign} p. 85-86 for the formulation), so
$$J_n\leq(1+o(1))\sqrt{2\pi}\frac{g_n(t_0)t_0}{\sqrt{n}}.$$
On the other hand, since $\varphi(t)$ is nondecreasing and positive,
\begin{equation}\label{left_int}
J_n\geq \int_0^{t_0} t^n e^{-\varphi(t)}dt\geq e^{-\varphi(t_0)}\int_0^{t_0} t^n dt= \frac{t_0 g_n(t_0)}{n+1}.
\end{equation}
$\square$

The next Lemma is a simple fact which we shall apply it to estimate the ``tails'' of $J_n$.
\begin{lemma}\label{log_concave}
Let $g(t)=e^{f(t)}$ be a log concave function on [a,b] (where both $a$ and $b$ may be infinite). Assume that $f\in C^2[a,b]$ and that $t_0$ is the unique point of maxima of $f(t)$. Assume that $t_0>0.$ Consider $x>0$ and $\psi>0$ such that
$$f(t_0)-f((1+x)t_0)\geq \psi.$$
Then,
$$\int_{(1+x)t_0}^{b} g(t)dt\leq \frac{x t_0 g(t_0)}{\psi e^{\psi}}.$$
Similarly, if $f(t_0)-f((1-x)t_0)\geq \psi,$
$$\int_a^{(1-x)t_0} g(t)dt\leq \frac{x t_0 g(t_0)}{\psi e^{\psi}}.$$
\end{lemma}
\noindent\textbf{Proof.} We pick any $t>(1+x)t_0$. First, we notice by concavity:
$$\psi\leq f(t_0)-f((1+x)t_0)\leq -f'((1+x)t_0)x t_0.$$
Next, since $f(t)$ is concave,
$$f(t)\leq f'((1+x)t_0)(t-(1+x)t_0)+f((1+x)t_0)\leq$$
$$-\frac{\psi}{x t_0}(t-(1+x)t_0)+f(t_0)-\psi.$$
Thus, for $t> (1+x)t_0$,
\begin{equation}\label{soslatsa}
g(t)\leq g(t_0) e^{-\psi}e^{{-\frac{\psi}{x t_0}}(t-(1+x)t_0)}.
\end{equation}
Consequently,
$$\int_{(1+x)t_0}^{b} g(t)dt\leq g(t_0) e^{-\psi}\int_0^{\infty} e^{{-\frac{\psi}{x t_0}}s}ds\leq \frac{x t_0 g(t_0)}{\psi e^{\psi}}.$$
The second part of the Lemma can be obtained similarly. $\square$

We note, that the condition $t_0>0$ in the above Lemma is not crucial, and everything can be restated for $t_0<0$. For our purposes it is enough to consider $t_0> 0.$

The function $g_n(t)=t^n e^{-\varphi(t)}$ is log concave on $[0,\infty)$, and we shall apply Lemma \ref{log_concave} with $g(t)=g_n(t)$ and $\psi=1$.
\begin{defn}\label{lambdas}
Define the ``outer'' $\lambda_o$ to be a positive number satisfying:
\begin{equation}\label{lambda_def_outer}
\varphi(t_0(1+\lambda_o))-\varphi(t_0)-n\log(1+\lambda_o)=1.
\end{equation}
Similarly, define the ``inner'' $\lambda_i$ as follows:
\begin{equation}\label{lambda_def_inner}
\varphi(t_0(1-\lambda_i))-\varphi(t_0)-n\log(1-\lambda_i)=1.
\end{equation}
We put
\begin{equation}\label{lambda}
\lambda:=\lambda_i+\lambda_o.
\end{equation}
\end{defn}

We note that (\ref{lambda_def_outer}) is equivalent to
\begin{equation}\label{lambda_def_outer1}
g_n(t_0)=e g_n(t_0(1+\lambda_o)),
\end{equation}
and (\ref{lambda_def_inner}) is equivalent to
\begin{equation}\label{lambda_def_inner1}
g_n(t_0)=e g_n(t_0(1-\lambda_i)).
\end{equation}

Parameter $\lambda$ from (\ref{lambda}) has a nice property:
\begin{lemma}\label{int_lam}
$$J_n \approx\lambda t_0 g_n(t_0).$$
\end{lemma}
\noindent\textbf{Proof.} We apply the first part of Lemma \ref{log_concave} with $x=\lambda_o$ and $\psi=1$. We get
\begin{equation}\label{1}
\int_{t_0(1+\lambda_o)}^{\infty} g_n(t)dt\leq \frac{1}{e}\lambda_o t_0 g_n(t_0).
\end{equation}
Similarly, the second part of the Lemma applied with $x=\lambda_i$, gives
\begin{equation}\label{2}
\int_0^{t_0(1-\lambda_i)} g_n(t)dt\leq \frac{1}{e}\lambda_i t_0 g_n(t_0).
\end{equation}
Along with the above, we observe:
\begin{equation}\label{3}
\int_{t_0(1-\lambda_i)}^{t_0(1+\lambda_o)} g_n(t)dt\leq (\lambda_i+\lambda_o) t_0 g_n(t_0)=\lambda t_0 g_n(t_0).
\end{equation}
From (\ref{1}), (\ref{2}) and (\ref{3}), applied together with the definition of $\lambda$, it follows that:
\begin{equation}\label{above1}
J_n\leq \frac{e+1}{e}\lambda t_0 g_n(t_0).
\end{equation}
On the other hand,
$$J_n\geq \int_{t_0(1-\lambda_i)}^{t_0(1+\lambda_o)} g_n(t)dt\geq$$
\begin{equation}\label{below1}
\lambda_i t_0 g_n((1+\lambda_i)t_0)+\lambda_o t_0 g_n((1+\lambda_o)t_0)=\frac{1}{e}\lambda t_0 g_n(t_0),
\end{equation}
where the last equality is obtained in a view of (\ref{lambda_def_inner1}) and (\ref{lambda_def_outer1}). $\square$

\begin{remark}\label{lambdy}
Let us note, that Lemma \ref{mainintegral} together with (\ref{above1}) and (\ref{below1}) leads to the estimates:
\begin{equation}\label{lambexact}
\frac{e}{e+1}\frac{1}{n+1}\leq \lambda \leq (1+o(1))\sqrt{2\pi}e\frac{1}{\sqrt{n}}.
\end{equation}

The above implies also, that both ``inner'' and ``outer'' lambdas are asymptotically bounded by $\frac{1}{\sqrt{n}}$.
In addition $\lambda_i\succsim \frac{1}{n}$. To see this, we write:
\begin{equation}\label{eq1}
\int_0^{t_0} t^n e^{-\varphi(t)}dt\geq e^{-\varphi(t_0)}\int_0^{t_0} t^n dt= \frac{t_0 g_n(t_0)}{n+1}.
\end{equation}
On the other hand, we estimate:
$$\int_{t_0(1-\lambda_i)}^{t_0}t^n e^{-\varphi(t)}dt\leq \lambda_i t_0 g_n(t_0).$$
Finally, we use (\ref{2}) and conclude:
$$\int_0^{t_0} t^n e^{-\varphi(t)}dt\leq \int_0^{t_0(1-\lambda_i)}t^n e^{-\varphi(t)}dt+\int_{t_0(1-\lambda_i)}^{t_0}t^n e^{-\varphi(t)}dt\leq $$
\begin{equation}\label{eq2}
\frac{1}{e}\lambda_i t_0 g_n(t_0)+\lambda_i t_0 g_n(t_0).
\end{equation}
The inequalities (\ref{eq1}) and (\ref{eq2}) yield the estimate $\lambda_i\geq \frac{e}{e+1}\frac{1}{n+1}$.

However, $\lambda_o$ can be arbitrarily small: for any $\varepsilon>0$ there exist a measure with continuous density (close to the one of the normalized Lebesgue measure on the unit ball) so that $\lambda_o<\varepsilon$.

\end{remark}

\begin{remark}\label{bigintegrals}
Remark \ref{lambdy} shows that $\lambda_o$ and $\lambda_i$ are $o(1)$ when $n\rightarrow \infty$. Consequently, for sufficiently large $n$,
$$\int_{\frac{1}{2} t_0}^{2 t_0} g_n(t)dt\approx J_n,$$
\end{remark}

The following fact is believed to be well known (see Remark 3.4 from \cite{bob2} for the best possible estimate).
\begin{lemma}\label{moments}
For all $n\geq 2$,
$$\frac{J_{n}}{J_{n-1}}\approx t_0.$$
\end{lemma}
\noindent\textbf{Proof.} In a view of Remark \ref{bigintegrals},
$$J_n=\int_0^{\infty} t^n e^{-\varphi(t)}dt\approx \int_{\frac{1}{2}t_0}^{2t_0} t^n e^{-\varphi(t)}dt\approx$$
$$t_0\int_{\frac{1}{2}t_0}^{2t_0} t^{n-1} e^{-\varphi(t)}dt\approx t_0 J_{n-1},$$
which completes the proof of the Lemma. $\square$

Let us consider some computable examples of $\gamma$-surface area. The first natural example to look at is the sphere of radius $R>0$.
$$\gamma(R\mathbb{S}^n)=\frac{1}{(n+1)\nu_{n+1}J_n}\int_{R\mathbb{S}^{n}} e^{-\varphi(|y|)}d\sigma(y)=\frac{|R\mathbb{S}^n|e^{-\varphi(R)}}{(n+1)\nu_{n+1}J_n}=$$
$$\frac{R^n e^{-\varphi(R)}}{J_n}\approx\frac{g_n(R)}{\lambda t_0 g_n(t_0)}.$$
Since $t_0$ is the maximum point for $g_n(t_0)$, we notice that among all the spheres, $t_0 \mathbb{S}^n$ has the maximal $\gamma$-surface area,
and it is equivalent to $\frac{1}{\lambda t_0}$.

Next, for a unit vector $\theta$ we consider the half space $H_{\theta}=\{y:\,\langle y,\theta\rangle\leq 0\}$.
\begin{equation}\label{hyp}
\gamma(\partial H_{\theta})=\frac{1}{(n+1)\nu_{n+1}J_n}\int_{\R^n} e^{-\varphi(|y|)} dy=\frac{n\nu_n J_{n-1}}{(n+1)\nu_{n+1}J_{n}}.
\end{equation}

Applying the fact that $\frac{\nu_n}{\nu_{n+1}}=\frac{\sqrt{n}}{\sqrt{2\pi}}(1+o(1))$ together with Lemma \ref{moments} and (\ref{hyp}), we obtain that $\gamma(H)\approx\frac{\sqrt{n}}{t_0}$.

We shall use a trick from \cite{ball} to show a rough upper bound for $\gamma(\partial Q)$.

\begin{lemma}\label{rough}
$\gamma(\partial Q)\precsim \frac{n}{t_0}$ for any convex set $Q$.
\end{lemma}
\noindent\textbf{Proof.}
We obtain the following integral expression for the density:
$$e^{-\varphi(|y|)}=\int_{|y|}^{\infty} \varphi'(t)e^{-\varphi(t)}dt=\int_0^{\infty} \varphi'(t)e^{-\varphi(t)}\chi_{[0,t]}(|y|)dt,$$
where $\chi_{[0,t]}$ stands for characteristic function of the interval $[0,t]$. Thus
$$\gamma(\partial Q)=\frac{1}{(n+1)\nu_{n+1} J_n}\int_{\partial Q}\int_0^{\infty} \varphi'(t)e^{-\varphi(t)}\chi_{[0,t]}(|y|)dt d\sigma(y)=$$
$$\frac{1}{(n+1)\nu_{n+1} J_n}\int_0^{\infty} \varphi'(t)e^{-\varphi(t)}|\partial Q\cap tB_2^{n+1}|dt,$$
which by can be estimated from above by
\begin{equation}\label{byparts}
\frac{(n+1)\nu_{n+1}}{(n+1)\nu_{n+1} J_n}\int_0^{\infty} t^n \varphi'(t)e^{-\varphi(t)}dt,
\end{equation}
since $ Q\cap tB_2^{n+1}\subset tB_2^{n+1}$, and thus $|\partial Q\cap tB_2^{n+1}|\leq |\partial tB_2^{n+1}|$ by convexity.

After integrating (\ref{byparts}) by parts and applying Lemma \ref{moments}, we get
$$\gamma(\partial Q)\leq n\frac{J_{n-1}}{J_n}\approx \frac{n}{t_0}. \square$$

The next lemma is an important tool in our proof.

\begin{lemma}\label{annulus}
Assume that there exists a positive $\mu$ such that
\begin{equation}\label{condition_lem_an}
\varphi\left(t_0(1+\mu)\right)-\varphi(t_0)-n\log(1+\mu)\geq\log\left({\mu \sqrt{\frac{n}{\lambda}}}\right)\geq 1.
\end{equation}
Define
$$A:=(1+\mu)t_0 B_2^{n+1}\setminus \frac{t_0}{e}B_2^{n+1}.$$
Then
$$\gamma(\partial Q\setminus A)\precsim \frac{\sqrt{n}}{t_0\sqrt{\lambda}}.$$
\end{lemma}
\noindent\textbf{Proof.} First, define the surface $B=\partial Q\cap \frac{t_0}{e}B_2^{n+1}$. Then,
$$\gamma(B)= \frac{1}{(n+1) \nu_{n+1} J_n} \int_{B} e^{-\varphi(|y|)} d\sigma(y)\leq \frac{|B|}{(n+1) \nu_{n+1} J_n}\leq$$
\begin{equation}\label{estimate1}
\frac{|\frac{t_0}{e} \mathbb{S}^n|}{(n+1) \nu_{n+1} J_n}\approx \frac{t_0^n}{e^n \lambda t_0 e^{-\varphi(t_0)}t_0^n}=\frac{1}{\lambda t_0} \frac{e^{\varphi(t_0)}}{e^n},
\end{equation}
where the equivalence follows from Lemma \ref{int_lam}. Recalling (\ref{philessthann}), which states that $\varphi(t_0)\leq n$, we estimate (\ref{estimate1})
from above by $\frac{1}{\lambda t_0}$. We recall as well, that $\frac{1}{\lambda t_0}\approx\gamma(t_0 \mathbb{S}^n)\precsim \frac{\sqrt{n}}{t_0\sqrt{\lambda}}$,
since $\lambda\succsim \frac{1}{n}$.

Next, let the surface $P=\partial Q\setminus (1+\mu)t_0 B_2^{n+1}.$ As in Lemma \ref{rough}, we make use of the estimate (\ref{byparts}) and integrate by parts:
$$\gamma(P)\leq \frac{1}{J_n}\int_{(1+\mu)t_0}^{\infty} t^n \varphi'(t)e^{-\varphi(t)}dt\precsim$$
\begin{equation}\label{mu_comp}
\frac{g_n((1+\mu)t_0)+n\int_{(1+\mu)t_0}^{\infty} g_{n-1}(t)dt}{\lambda t_0 g_n(t_0)}.
\end{equation}
Lemma \ref{log_concave}, applied with $x=\mu$ and $\psi=\log\left(\mu\sqrt{\frac{n}{\lambda}}\right)$, entails that (\ref{mu_comp}) is less than
$$
\frac{e^{-\psi}}{\lambda t_0}+\frac{n\mu}{\lambda t_0\psi e^{\psi}}= \frac{1}{\lambda t_0} \times (1+\frac{\mu n}{\psi}) e^{-\psi}\precsim \frac{\sqrt{n}}{t_0\sqrt{\lambda}},
$$
where the last bound follows if we plug $\psi=\log\left(\mu\sqrt{\frac{n}{\lambda}}\right)$ and use the fact that $\psi\geq 1$. We also used Remark \ref{lambdy} which yields the fact that $\frac{1}{\lambda t_0}\precsim \frac{\sqrt{n}}{t_0\sqrt{\lambda}}.$ $\square$

The next Lemma shows, that $\mu$ in Lemma \ref{annulus} can be chosen very small.

\begin{lemma}\label{mu}
$$\mu=\frac{\log n}{\sqrt{n}}$$
satisfies the condition of Lemma \ref{annulus} for sufficiently large $n$.
\end{lemma}
\noindent\textbf{Proof.} First, notice that $\varphi((1+\mu)t_0)-\varphi(t_0)\geq \varphi'(t_0)\mu t_0=n\mu.$ Thus
\begin{equation}\label{dot}
\varphi(t_0(1+\mu))-\varphi(t_0)-n\log(1+\mu)\geq n(\mu-\log(1+\mu)).
\end{equation}
Plugging $\mu=\frac{\log n}{\sqrt{n}}$ into (\ref{dot}) and applying the Taylor approximation for logarithm, we get that the right hand side of (\ref{dot}) is approximately equal to
\begin{equation}\label{right}
\sqrt{n}\log{n}-n\log\left(1+\frac{\log n}{\sqrt{n}}\right)=\frac{\log^2 n}{2}+o(1).
\end{equation}
In order to satisfy (\ref{condition_lem_an}), we need to estimate $\log \left(\mu\sqrt{\frac{n}{\lambda}}\right)$ from above:
\begin{equation}\label{left}
\log \left(\mu\sqrt{\frac{n}{\lambda}}\right)=\log\left(\frac{\log n}{\sqrt{n}}\sqrt{\frac{n}{\lambda}}\right)\leq \log(5n\log n),
\end{equation}
since $\lambda\geq \frac{e}{e+1}\frac{1}{n}$ (see Remark \ref{lambdas}).
Observing, that for all $n\geq 12$, $\log\left(5n\log n\right)\leq \frac{\log^2 n}{2}+o(1)$, we obtain the Lemma. $\square$

\section {Connections to Probability}
We consider a random vector $X$ in $\R^{n+1}$ distributed with respect to $\gamma$. Then $|X|$ is a random variable distributed on $[0,\infty)$ with density $\frac{g_n(t)}{J_n}$. We shall use standard notation for its expectation and variance: $E=\mathbb{E}|X|=\frac{1}{J_n}\int_0^{\infty} t g_n(t)dt$ and
\begin{equation}\label{var}
\sigma^2=Var|X|=\frac{1}{J_n} \int_0^{\infty} (t-E)^2 g_n(t)dt.
\end{equation}
The next two Lemmas give an expression for the expectation and variance of $|X|$ in terms of our parameters $\lambda$ and $t_0$, which will be used to restate Theorem \ref{th:main}.
\begin{lemma}\label{expectationnn}
$$\mathbb{E}|X|\approx t_0.$$
\end{lemma}
\noindent\textbf{Proof.} We write
$$\mathbb{E}|X|=[(n+1)\nu_{n+1} J_{n}]^{-1}\int_{\R^{n+1}} |y|e^{-\varphi(|y|)}dy=$$
$$\frac{1}{J_n} \int_0^{\infty} t^{n+1} e^{-\varphi(t)}dt=\frac{J_{n+1}}{J_n}\approx t_0,$$
where the last equivalence follows from Lemma \ref{moments}. $\square$

\begin{lemma}\label{varience}
$$Var|X|\approx (\lambda t_0)^2.$$
\end{lemma}
\noindent\textbf{Proof.} We notice first that (\ref{var}) implies:
\begin{equation}\label{11111}
\int_0^{\infty} g_n(t)\frac{(t-E)^2}{4\sigma^2}dt=\frac{J_n}{4}.
\end{equation}
Subtracting (\ref{11111}) from the equation $J_n=\int_0^{\infty} g_n(t)dt$, we get
$$\int_0^{\infty} g_n(t)\left(1-\frac{(t-E)^2}{4\sigma^2}\right)dt=\frac{3}{4} J_n.$$
We observe that $1-\frac{(t-E)^2}{4\sigma^2}$ is between zero and one whenever $|t-E|\leq 2\sigma$, and negative otherwise. Thus
\begin{equation}\label{above}
\int_{E-2\sigma}^{E+2\sigma} g_n(t)dt\geq \int_{E-2\sigma}^{E+2\sigma} g_n(t)\left(1-\frac{(t-E)^2}{4\sigma^2}\right)dt\geq \frac{3}{4} J_n.
\end{equation}
On the other hand,
$$\int_{E-2\sigma}^{E+2\sigma} g_n(t)dt\leq 4\sigma\times \max_{t\in [E-2\sigma,E+2\sigma]} g_n(t)\leq$$
\begin{equation}\label{below}
 4\sigma \max_{t\in [0,\infty)} g_n(t)=4\sigma g_n(t_0).
\end{equation}
Bringing together Lemma \ref{int_lam}, (\ref{above}) and (\ref{below}), we get
$$4\sigma g_n(t_0)\geq \frac{3}{4} J_n\approx \lambda t_0 g_n(t_0),$$
and thus $\sigma\succsim \lambda t_0.$

Next, we shall obtain the reverse estimate. We note that the expression
$$\int_{0}^{\infty} (t-\tau)^2 g_n(t)dt$$
is minimal when $\tau=E$. Thus for $\tau=t_0(1+\lambda)$ we get:
$$\sigma^2J_n\leq \int_{0}^{\infty} (t-t_0(1+\lambda))^2 g_n(t)dt=$$
\begin{equation}\label{split}
\int_{0}^{t_0(1-\lambda)}+\int_{t_0(1-\lambda)}^{t_0(1+\lambda)}+\int_{t_0(1+\lambda)}^{\infty} (t-t_0(1+\lambda))^2 g_n(t)dt.
\end{equation}
The second integral in (\ref{split}) can be bounded by
\begin{equation}\label{second}
\max_{t\in [t_0-\lambda t_0,t_0+\lambda t_0]} (t-t_0(1+\lambda))^2\int_{t_0(1-\lambda)}^{t_0(1+\lambda)} g_n(t)dt\precsim (\lambda t_0)^2 J_n.
\end{equation}

In order to estimate the third integral we apply (\ref{soslatsa}) with $g(t)=g_n(t)$, $\psi= 1$ and $x=\lambda$. It implies that for all $t>t_0(1+\lambda),$ the following holds:
$$g_n(t)\precsim g_n(t_0) e^{-\frac{1}{\lambda t_0}(t-t_0(1+\lambda))}.$$
Thus the third integral from (\ref{split}) can be estimated from above with
$$g_n(t_0)\int_{t_0(1+\lambda)}^{\infty} (t-t_0(1+\lambda))^2 e^{-\frac{1}{\lambda t_0}(t-t_0(1+\lambda))}dt=$$
$$(\lambda t_0)^3 g_n(t_0)\int_0^{\infty} s^2 e^{-s} ds=2(\lambda t_0)^2\lambda t_0 g_n(t_0)\approx (\lambda t_0)^2 J_n,$$
where the last equivalence follows from Lemma \ref{int_lam}. The first integral in (\ref{split}) can be estimated similarly (with the loss of $e^{-2}$). Adding both of them together with (\ref{second}), we obtain that
$$\sigma^2 J_n\precsim (\lambda t_0)^2 J_n,$$
which finishes the proof. $\square$

Now we are ready to restate Theorem \ref{th:main}:

\begin{theorem}\label{main}
Fix $n\geq 2$. Let $t_0$ be the solution of $\varphi'(t)t=n-1.$
Define $\widetilde{\lambda}=\frac{\int_0^{\infty} t^{n-1} e^{-\varphi(t)}dt}{t_0^n e^{-\varphi(t_0)}}.$ Then
$$\max_{Q\in \mathcal{K}_n}\gamma(\partial Q)\approx\frac{\sqrt{n}}{\sqrt{\widetilde{\lambda}} t_0}.$$
\end{theorem}

From now on we will be after proving Theorem \ref{main}. Notice, that by Lemma \ref{int_lam}, $\widetilde{\lambda}$ is equivalent to $\lambda$, defined in the previous section.

\begin{remark}
The statement of Theorem \ref{main} becomes shorter if the measure is isotropic. We refer to \cite{MP} and \cite{Kl} for the definitions and details.
Here we observe only, that $t_0=\sqrt{n}$ for isotropic measures on $\R^n$, and after making a change of variables $\widetilde{\varphi}(t)=\varphi(\frac{t_0}{\sqrt{n}}t)$, we get a measure $\widetilde{\gamma}$ with density $C(n)e^{-\widetilde{\varphi}(|y|)}$, which has properties similar to $\gamma$ and for which the statement of Theorem \ref{main} becomes:
$$\max_{Q\in \mathcal{K}_n}\widetilde{\gamma}(\partial Q)\approx \frac{1}{\sqrt{\widetilde{\lambda}}}.$$
\end{remark}

\begin{remark}
For $p\geq 1$ we define $\gamma_p$ to be a probability measure on $\R^n$ with density $C_{n,p}e^{-\frac{|y|^p}{p}}$ (as in (\ref{ya}). In this case $\varphi(t)=\frac{t^p}{p}$,
and $\varphi'(t)t=t^p$. Thus, for such measures $t_0=(n-1)^{\frac{1}{p}}$ (see (\ref{todef}) for the definition of $t_0$). Also, Laplace method entails, that
$$J_n=c(p) \frac{(n-1)^{\frac{n}{p}} e^{-\frac{n-1}{p}}}{\sqrt{n}}=c(p)\frac{t_0 g_n(t_0)}{\sqrt{n}}.$$
(see \cite{L} for the details.) In a view of Lemma \ref{int_lam} we conclude, that in this case $\lambda\approx \frac{1}{\sqrt{n}}$. So Theorem \ref{main} asserts, that
$$\max_{Q\in \mathcal{K}_n} \gamma_p(\partial Q) \approx C(p)n^{\frac{3}{4}-\frac{1}{p}},$$
which means that the result of \cite{L} for the case $p\geq 1$,  the result of \cite{fedia} for the standard Gaussian measure, and the result from \cite{ball} are consequences of the current one.
\end{remark}

\section {Upper bound}

We will use the approach developed by Nazarov in \cite{fedia}. We pick a convex set $Q$. The aim is to estimate $\gamma(\partial Q)$ from above. By log concavity of measure $\gamma$, we may assume that $Q$ contains the origin: otherwise we may shift $Q$ towards the origin so that the surface area does not decrease. Indeed, if $Q$ does not contain the origin, let $y_0\in Q$ be the closest point to the origin. Apply the shift $S(y)=y-y_0.$ The body $S(Q)$ contains the origin in it's boundary, and also $|y-y_0|\leq |y|$ for all $y\in Q$. Since $\varphi$ is increasing, we get $\varphi(|y-y_0|)\leq\varphi(|y|)$, and thus $\gamma (\partial S(Q))\geq \gamma(\partial Q)$. Moreover, by continuity of $\varphi(t)$ we may assume that the origin is contained not in the boundary, but in the interior of $Q$.

Let us consider ``polar'' coordinate system $x=X(y,t)$ in $\R^{n+1}$ with $y\in \partial Q$, $t>0$. We write
$$C_{n+1}\int_{\R^n} \F d\sigma(y)=C_{n+1}\int_0^{\infty}\int_{\partial Q} D(y,t)e^{-\varphi(|X(y,t)|)} d\sigma(y)dt,$$
where $D(y,t)$ is the Jacobian of $x\rightarrow X(y,t)$. Define
\begin{equation}\label{xi}
\xi(y)=e^{\varphi(|y|)}\int_0^{\infty} D(y,t)e^{-\varphi(|X(y,t)|)} dt.
\end{equation}
Then
$$1=C_{n+1}\int_{\partial Q} \F\xi(y)d\sigma(y),$$
and thus
\begin{equation}\label{upperbound}
\gamma(\partial Q)=C_{n+1}\int_{\partial{Q}}\F dy\leq \frac{1}{\min\limits_{y\in\partial Q} \xi(y)}.
\end{equation}
Following \cite{fedia}, we shall consider two such systems.

\subsection{First coordinate system}

We consider ``radial'' polar coordinate system $X_1(y,t)=yt$. The Jacobian $D_1(y,t)=t^{n}|y|\alpha$, where
\begin{equation}\label{alpha}
\alpha=\alpha(y)=\cos(y,\n_y),
\end{equation}
where $\n_y$ stands for a normal vector at $y$. Without loss of generality we assume that $n_y$ is defined uniquely for every $y\in \partial Q$. Rewriting (\ref{xi}), making a change of variables $\tau=t|y|$ and applying Lemma \ref{moments}, we get:
$$\xi_1(y):=e^{\varphi(|y|)}\int_0^{\infty} t^{n}|y|\alpha e^{-\varphi(|ty|)} dt=$$
\begin{equation}\label{xi11}
e^{\varphi(|y|)}\alpha |y|^{-n}J_{n}\succsim t_0\alpha\lambda \frac{g_n(t_0)}{g_n(|y|)}.
\end{equation}

We define $x=x(y)$ to satisfy $|y|=(1+x)t_0$ and
\begin{equation}\label{psi}
\psi(x):=\varphi((1+x)t_0)-\varphi(t_0)-n\log(1+x)=\log \frac{g_n(t_0)}{g_n((1+x)t_0)}.
\end{equation}
Then, by (\ref{xi11}),
\begin{equation}\label{xi1}
\xi_1(y)\succsim t_0\alpha\lambda e^{\psi(x)}.
\end{equation}

\begin{remark}
For the sake of completeness we note, that the above formula might as well be obtained by projecting the set on the unit sphere and passing to new coordinates. Indeed, let $x=\frac{y}{|y|}$. Then the coordinate change writes as $d\sigma(y)=\frac{|y|^n}{\alpha(y)}d\sigma(x)$, and we obtain
$$\gamma(\partial Q)=\C\int_{\partial Q} \F d\sigma(y)=$$
$$\C\int_{\mathbb{S}^n} \F \frac{|y|^n}{\alpha(y)}d\sigma(x)\leq \max_{y\in\partial Q}{\frac{g_n(|y|)}{\alpha(y)J_n}},$$
which is equivalent to the bound we obtain from (\ref{xi}) and (\ref{xi1}). This observation shows, that no volume argument of the type (\ref{upperbound}) is needed
 here. However, we shall need it below.
\end{remark}

\subsection{Second coordinate system}
We consider ``normal'' polar coordinate system $X_2(y,t)=y+t\n_y$. Then $D_2(y,t)\geq 1$ for all $y\not\in Q$. We write
$$\varphi(|X_2(y,t)|)=\varphi(|y+t\n_y|)=\varphi\left(\sqrt{|y|^2+t^2+2t|y|\alpha}\right),$$
where $\alpha=\alpha(y)$ was defined by (\ref{alpha}). Let $\xi_2(y)$ be $\xi(y)$ from (\ref{xi}), corresponding to $X(y,t)=X_2(y,t)$. Then
\begin{equation}\label{startxi2}
\xi_2(y)\geq e^{\varphi(|y|)}\int_0^{\infty}e^{-\varphi\left(\sqrt{|y|^2+t^2+2t|y|\alpha}\right)}dt.
\end{equation}
Define $t_1=t_1(y)$ to be the largest number such that:
$$\varphi\left(\sqrt{|y|^2+t_1^2+2t_1|y|\alpha}\right)-\varphi(|y|)=1.$$
Such number always exists, since the function $\varphi\left(\sqrt{|y|^2+t^2+2t|y|\alpha}\right)$ is a nondecreasing continuous function of $t$ on $[0,\infty)$, and
$$\lim_{t\rightarrow +\infty} \varphi\left(\sqrt{|y|^2+t^2+2t|y|\alpha}\right)=+\infty.$$
We shall use an elementary inequality
$$\int f(x) d\mu(x)\geq a\times \mu(f(x)\geq a),$$
which holds for all positive integrable functions $f$.
Notice, that
$$|\{t\geq 0:\,e^{-\varphi\left(\sqrt{|y|^2+t^2+2t|y|\alpha}\right)}\geq e^{-\varphi(|y|)-1} \}|=t_1.$$
Thus the right hand side of (\ref{startxi2}) is asymptotically bounded from below by $t_1$.


We define $\Lambda(t):\, [0,\infty)\rightarrow [0,\infty)$ the relation
\begin{equation}\label{Lambdabig}
\varphi((1+\Lambda(t))t)-\varphi(t)=1.
\end{equation}

By the definition of $t_1=t_1(y)$,
$$\sqrt{1+\frac{t_1^2}{|y|^2}+\frac{2t_1\alpha}{|y|}}=1+\Lambda(|y|).$$
We solve the quadratic equation and obtain, that for all $y\in \partial Q$
\begin{equation}\label{xi2}
\xi_2(y)\geq \frac{t_1}{e} \succsim \frac{|y|\sqrt{\Lambda(|y|)+\Lambda^2(|y|)}}{\frac{\alpha(y)}{\sqrt{\Lambda(|y|)+\Lambda^2(|y|)}}+1}.
\end{equation}

\subsection{Cases.}
We shall split the space into several annuli and estimate $\gamma$-surface area of $\partial Q$ intersected with each annulus separately. The proof splits into several cases. Below we assume that $y\in\partial Q.$
\\
\\
\textbf{Case 1: $|y|\leq \frac{1}{2e}t_0$ or $|y|\geq (1+\frac{\log n}{\sqrt{n}})t_0$.}
\\
We define $\partial Q_1=\{y\in\partial Q: |y|\leq \frac{1}{2e}t_0\,\,\, or\,\,\, |y|\geq (1+\frac{\log n}{\sqrt{n}})t_0\}$.
Direct application of Lemmas \ref{annulus} and \ref{mu} asserts that the desired upper bound holds for $\gamma(\partial Q_1)$ (we remark, that even though the application of Lemma \ref{mu} requires $n\geq 12,$ we may apply Lemma \ref{rough} for $n\leq 12$ and select the proper constant at the end).
\\
\\
\textbf{Case 2: $\frac{1}{2e}t_0\leq |y|\leq (1-\frac{1}{n})t_0$.}\\
We define $\partial Q_2=\{y\in\partial Q: \frac{1}{2e}t_0\leq |y|\leq (1-\frac{1}{n})t_0\}$. Pick $y\in \partial Q_2.$ We observe:
$$\varphi\left((1-\frac{1}{n})t_0\times(1+\frac{1}{n})\right)-\varphi\left((1-\frac{1}{n})t_0\right)\leq $$
\begin{equation}\label{mvt}
\varphi(t_0)-\varphi\left((1-\frac{1}{n})t_0\right)\leq \frac{t_0}{n}\varphi'(t_0)=1.
\end{equation}
This asserts that $\Lambda((1-\frac{1}{n})t_0)\geq\frac{1}{n}$. We note, that $\Lambda(t)$ decreases, when $t$ increases. Thus $\Lambda(|y|)\geq \frac{1}{n}$ for any $y$ such that $|y|\leq (1-\frac{1}{n})t_0$. We rewrite (\ref{xi2}) and get the estimate
\begin{equation}\label{xi2case2start}
\xi_2(y) \succsim \frac{|y|}{\sqrt{n}}\times\frac{1}{\alpha\sqrt{n}+1}.
\end{equation}
Since $|y|$ is assumed to be asymptotically equivalent to $t_0$, (\ref{xi2case2start}) rewrites as
\begin{equation}\label{xi2case2}
\xi_2(y) \succsim \frac{t_0}{\sqrt{n}}\times\frac{1}{\alpha\sqrt{n}+1}.
\end{equation}
As for the first system, we apply a rough estimate $\psi(x)\geq 0$ and rewrite (\ref{xi1}) as follows:
\begin{equation}\label{xi1case2}
\xi_1(y)\succsim t_0\alpha\lambda.
\end{equation}
 We consider
\begin{equation}\label{xidef}
\xi(y):=\xi_1(y)+\xi_2(y)\succsim
\end{equation}
$$t_0\alpha\lambda+\frac{t_0}{\sqrt{n}}\times\frac{1}{\alpha\sqrt{n}+1}.$$
We minimize the above expression with respect to $\alpha\in[0,1]$.
The minimum is attained when $\alpha=\frac{1}{\sqrt{\lambda n}}$, and thus
$$\xi(y)\succsim \frac{t_0 \sqrt{\lambda}}{\sqrt{n}},$$
which together with (\ref{xi}) and (\ref{upperbound}) leads to the desired estimate for $\gamma(\partial Q_2)$.
\\
\\
\textbf{Case 3: $(1-\frac{1}{n})t_0\leq |y|\leq t_0$}\\
We define $\partial Q_3=\{y\in\partial Q: (1-\frac{1}{n})t_0\leq |y|\leq t_0\}$.
Along the annulus the value of $\varphi(t)$ doesn't change that much. Namely, since $\varphi(t)$ is nondecreasing and by (\ref{mvt}),
$$\varphi\left(t_0(1-\frac{1}{n})\right)\in[\varphi(t_0)-1,\varphi(t_0)].$$
So for all $y\in\partial Q_3$, $\varphi(|y|)\approx\varphi(t_0)$. Thus we write
$$\gamma(\partial Q_3)=[(n+1)\nu_{n+1} J_{n}]^{-1} \int_{\partial Q_3} e^{-\varphi(|y|)}d\sigma(y)\approx$$
$$\frac{e^{-\varphi(t_0)}}{(n+1)\nu_{n+1}J_n} \int_{\partial Q_3} d\sigma(y)=\frac{e^{-\varphi(t_0)}}{(n+1)\nu_{n+1}J_n}|\partial Q_3|.$$
Since $Q_3$ is a convex body contained in $t_0 B_2^{n+1}$, we get $|\partial Q_3|\leq |t_0\mathbb{S}^n|$, so the above is less than
$$\frac{e^{-\varphi(t_0)|t_0\mathbb{S}^{n}|}}{(n+1)\nu_{n+1}J_n}=\frac{e^{-\varphi(t_0)}t_0^n }{J_n}\approx \frac{1}{\lambda t_0},$$
where the last equivalence is a direct application of Lemma \ref{int_lam}.
We conclude that the portion of any convex set in a very thin annulus around the maximal sphere is at least as small as the maximal sphere itself,
and, in particular, smaller than our desired upper bound.\\
\\
\textbf{Case 4: $t_0\leq |y|\leq (1+\frac{\log n}{\sqrt{n}})t_0$.}
\\
This case is the hardest one. We face the problem of controlling $\Lambda(y)$: there is no way to get a proper lower bound for it unless we ``step inside''
the set a little bit. Fortunately, Lemma \ref{homothety} shows that stepping not too far does not change $\gamma-$surface area too much. So we will be estimating
$\xi_2\left(\frac{|y|}{(1+\frac{1}{n})^2}\right)$ from below, rather than $\xi_2(y)$. The key estimate in all our computation is the following Proposition.

\begin{proposition}\label{Lambdaprop}
For any $y$ such that $|y|\in[t_0, (1+\frac{\log n}{\sqrt{n}})t_0]$,
$$\Lambda\left(\frac{|y|}{(1+\frac{1}{n})^2}\right) \succsim \frac{1}{n}\times \frac{1}{\psi(x)+1+o(1)},$$
where $\psi(x)$ is defined by (\ref{psi}), $\Lambda(t)$ is defined by (\ref{Lambdabig}) and $|y|=(1+x)t_0.$
\end{proposition}
\noindent\textbf{Proof.} We fix $|y|=(1+x)t_0$. The parameter $x$ in this case ranges
between $0$ and $\frac{\log n}{\sqrt{n}}$. Notice that by the Mean Value Theorem,
\begin{equation}\label{Lambda1}
\Lambda(|y|)\succsim \frac{1}{|y|\varphi'((1+\Lambda(|y|))|y|)}.
\end{equation}

For any $y$ such that $|y|\geq t_0$,
\begin{equation}\label{llc}
\Lambda(|y|)\leq\frac{\varphi((1+\Lambda(|y|))|y|)-\varphi(|y|)}{|y|\varphi'(|y|)}=\frac{1}{|y|\varphi'(|y|)}\leq \frac{1}{t_0\varphi'(t_0)}=\frac{1}{n}.
\end{equation}
Since $\varphi'(t)$ is nondecreasing, (\ref{Lambda1}) is greater than $\frac{1}{|y|\varphi'((1+\frac{1}{n})|y|)}.$ We apply (\ref{Lambda1}) with $|y|=\frac{1+x}{(1+\frac{1}{n})^2}t_0$:

\begin{equation}\label{lambda_phi_prime}
\Lambda\left(\frac{1+x}{(1+\frac{1}{n})^2}t_0\right)\succsim \frac{1}{(1+x)t_0\varphi'(\frac{(1+x)t_0}{1+\frac{1}{n}})}\approx \frac{1}{t_0\varphi'\left(\frac{(1+x)t_0}{1+\frac{1}{n}}\right)},
\end{equation}
where the last equivalence holds in the current range of $x$. Next, we write that
\begin{equation}\label{42}
\varphi'\left(\frac{(1+x)t_0}{1+\frac{1}{n}}\right)\leq \frac{\varphi\left((1+x)t_0\right)-\varphi\left(\frac{(1+x)t_0}{1+\frac{1}{n}}\right)}{(1+x)t_0-\frac{(1+x)t_0}{1+\frac{1}{n}}}.
\end{equation}
We note, that
\begin{equation}\label{computationwithx}
(1+x)t_0-\frac{(1+x)t_0}{1+\frac{1}{n}}=\frac{(1+x)t_0}{n+1}\approx \frac{t_0}{n}
\end{equation}
in the current range of $x$. We shall invoke the function $\psi(x)$. Applying its definition (\ref{psi}) in the numerator and (\ref{computationwithx}) in the denominator of (\ref{42}), we get that (\ref{42}) is equivalent to
\begin{equation}\label{refernow}
\frac{\psi(x)+n\log(1+x)+\varphi(t_0)-\varphi\left(\frac{(1+x)t_0}{1+\frac{1}{n}}\right)}{\frac{t_0}{n}}.
\end{equation}
Notice now, that by the Mean Value Theorem,
\begin{equation}\label{uzhass}
\varphi\left(\frac{(1+x)t_0}{1+\frac{1}{n}}\right)-\varphi(t_0)\succsim \varphi'(t_0)t_0 \left(x-\frac{1+o(1)}{n}\right)=nx-1+o(1).
\end{equation}
By (\ref{refernow}) and (\ref{uzhass}),
$$\varphi'\left(\frac{(1+x)t_0}{1+\frac{1}{n}}\right)\precsim \frac{n}{t_0}\times\left(\psi(x)+n\log(1+x)-nx+1+o(1)\right).$$
An elementary inequality $x\geq \log(1+x)$ entails that
\begin{equation}\label{aaa}
\varphi'\left(\frac{(1+x)t_0}{1+\frac{1}{n}}\right)\precsim\frac{n}{t_0}\left(\psi(x)+1+o(1)\right).
\end{equation}
Finally, by (\ref{aaa}) and (\ref{lambda_phi_prime}) we conclude
$$
\Lambda\left(\frac{1+x}{(1+\frac{1}{n})^2}t_0\right) \succsim \frac{1}{n}\times \frac{1}{\psi(x)+1+o(1)}.\square
$$


In the next few lines we use notation $\Lambda=\Lambda(\frac{|y|}{(1+\frac{1}{n})^2})$ for clarity of the presentation. We consider
$$\widetilde{\xi(y)}:=\xi_1(y)+\xi_2(\frac{y}{(1+\frac{1}{n})^2})\succsim $$
\begin{equation}\label{ono}
t_0\alpha\lambda e^{\psi(x)}+\frac{t_0\sqrt{\Lambda+\Lambda^2}}{2\frac{\alpha}{\sqrt{\Lambda+\Lambda^2}}+1}.
\end{equation}
First, we shall minimize (\ref{ono}) with respect to $\alpha$. It is minimized whenever
$$\alpha\approx\alpha_{min}:=\sqrt{\Lambda+\Lambda^2}\left(\frac{1}{\sqrt{e^{\psi(x)}\lambda}}-1\right).$$
Since $\psi$ is increasing on $(t_0,\infty)$, and due to our restrictions of the case 4, we may assume that
$$\psi(x)\leq \psi(t_0(1+\frac{\log n}{\sqrt{n}}))= \log \left(\sqrt{\frac{n}{\lambda}}x\right) \leq$$
$$\log \left(\frac{\log n}{\sqrt{n}}\sqrt{\frac{n}{\lambda}}\right)=\log\left( \frac{\log n}{\sqrt{\lambda}}\right).$$
Consequently,
$$\sqrt{\lambda e^{\psi(x)}}\leq \sqrt[4]{\lambda} \sqrt{\log n}=o(1),$$
and thus $\alpha_{min}\approx \sqrt{\frac{\Lambda+\Lambda^2}{e^{\psi(x)}\lambda}}$. Plugging it into (\ref{ono}), we obtain:
\begin{equation}\label{xi_first}
\widetilde{\xi(y)}\succsim t_0 \sqrt{\lambda (\Lambda+\Lambda^2) e^{\psi(x)}}.
\end{equation}

Finally, we apply (\ref{xi_first}) together with Proposition \ref{Lambdaprop}:
\begin{equation}\label{!}
\widetilde{\xi(y)}\succsim \frac{t_0 \sqrt{\lambda}}{\sqrt{n}}\sqrt{\frac{e^{\psi(x)}}{\psi(x)+1+o(1)}}\succsim \frac{t_0 \sqrt{\lambda}}{\sqrt{n}},
\end{equation}
where the last inequality holds since $\psi(x)$ is positive.
\subsection{Balancing for Case 4}
We restrict our attention on the part of the boundary which satisfies the condition of the Case 4. Namely, denote
$\partial Q_4:=\{y\in \partial Q: \, t_0\leq |y|\leq (1+\frac{\log n}{\sqrt{n}})t_0 \}$.

We would like to apply (\ref{xi}) and (\ref{upperbound}) with $\xi(y)=\xi_1(y)+\xi_2(y)$ and finish the proof, but unfortunately
we only have a lower bound for $\widetilde{\xi(y)}=\xi_1(y)+\xi_2(\frac{y}{(1+\frac{1}{n})^2})$. So we have to be a little bit more careful.
We define $A=\{y\in \partial Q_4:\xi_1(y)\geq\xi_2(\frac{y}{(1+\frac{1}{n})^2})\}$ and its compliment
$B=\{y\in \partial Q_4:\xi_1(y)<\xi_2(\frac{y}{(1+\frac{1}{n})^2})\}.$ Note, that both $A$ and $B$ are $\gamma-$measurable,
since $\xi_1$ and $\xi_2$ are Borell functions and $\gamma$ is absolutely continuous with respect to Lebesgue measure.
We shall apply (\ref{xi}) and (\ref{upperbound}) with $\xi(y)=\xi_1(y)$ on the set $A$ and with $\xi(y)=\xi_2(y)$ on the set $\frac{1}{(1+\frac{1}{n})^2}B$.

We write that
$$1\geq [(n+1)\nu_{n+1} J_{n}]^{-1}\int_{A}\int_0^{\infty} e^{-\varphi(X_1(y,t))} D_1(y,t) dt d\sigma(y)=$$
$$[(n+1)\nu_{n+1} J_{n}]^{-1}\int_{A}e^{-\varphi(|y|)} \xi_1(y) d\sigma(y)\geq \gamma(A) \min_{y\in A} \xi_1(y).$$
Thus,
\begin{equation}\label{111}
\gamma(A)\leq \frac{1}{\min_{y\in A} \xi_1(y)}.
\end{equation}
Similarly, we write
$$1\geq [(n+1)\nu_{n+1} J_{n}]^{-1} \int_{\frac{1}{(1+\frac{1}{n})^2}B}\int_0^{\infty} e^{-\varphi(X_2(y,t))} D_2(y,t) dt d\sigma(y)=$$
$$
[(n+1)\nu_{n+1} J_{n}]^{-1}\int_{\frac{1}{(1+\frac{1}{n})^2}B}e^{-\varphi(|y|)} \xi_2(y) d\sigma(y)\geq
$$
\begin{equation}\label{two}
\min_{y\in B} \xi_2\left(\frac{y}{(1+\frac{1}{n})^2}\right)\gamma\left(\frac{1}{(1+\frac{1}{n})^2}B\right).
\end{equation}
We apply Lemma \ref{homothety} for $M=\frac{1}{(1+\frac{1}{n})^2}B$ together with (\ref{two}), and conclude that
\begin{equation}\label{222}
\gamma(B)\precsim \frac{1}{\min_{y\in B} \xi_2(\frac{y}{(1+\frac{1}{n})^2})}.
\end{equation}
From (\ref{111}) and (\ref{222}) we obtain the following:
$$\gamma(\partial Q_4)=\gamma(A\cup B)\precsim\frac{1}{\min_{y\in A} \xi_1(y)}+\frac{1}{\min_{y\in B} \xi_2(\frac{y}{(1+\frac{1}{n})^2})}.$$
Invoking the definitions of the sets $A$ and $B$, we notice, that
$$\min_{y\in A} \xi_1(y)\geq \frac{1}{2} \min_{y\in A} \left(\xi_1(y)+\xi_2\left(\frac{y}{(1+\frac{1}{n})^2}\right)\right)\geq$$
$$\frac{1}{2} \min_{y\in \partial Q_4} \left(\xi_1(y)+\xi_2\left(\frac{y}{(1+\frac{1}{n})^2}\right)\right),$$
as well as
$$\min_{y\in B} \xi_2\left(\frac{y}{(1+\frac{1}{n})^2}\right)\geq \frac{1}{2} \min_{y\in B} \left(\xi_1(y)+\xi_2\left(\frac{y}{(1+\frac{1}{n})^2}\right)\right)\geq$$
$$\frac{1}{2} \min_{y\in \partial Q_4} \left(\xi_1(y)+\xi_2\left(\frac{y}{(1+\frac{1}{n})^2}\right)\right),$$
since the minimum over the smaller set is greater than the minimum over the larger set. We conclude, that
$$\gamma(\partial Q_4)\precsim \frac{1}{\min_{y\in \partial Q_4} \widetilde{\xi(y)}},$$
where $\widetilde{\xi(y)}=\xi_1(y)+\xi_2(\frac{y}{(1+\frac{1}{n})^2})$. The desired lower bound for this quantity was obtained earlier (\ref{!}),
which finishes the proof of the upper bound part for Theorem \ref{main}.

\section{Lower bound}
It seems impossible to construct an explicit example of a convex set $Q$ with $\gamma(\partial Q)\approx \frac{\sqrt{n}}{\sqrt{\lambda} t_0}$.
So we provide a probabilistic construction similar to the one in \cite{fedia}. Namely, we shall consider a random polytope circumscribed around a
sphere of a certain radius. The radius of the sphere and the number of faces shall be chosen so that most of the time $\alpha(y)=\cos(y,\n_y)\approx \alpha_{min}$
which appears in the proof of the upper bound, and so that large enough portion of the polytope falls close to the maximal sphere $t_0\mathbb{S}^n$.
As it was shown in Lemma \ref{int_lam}, a lot of the measure is concentrated in the thin annulus around $t_0\mathbb{S}^n$; more precise results
describing the decay outside of the annulus were obtained in \cite{Kl} (Theorem 1.4) and \cite{Kl2} (Theorem 4.4).
For simplicity of the calculations, we only look at the portion of the polytope in that annulus, and it turns out to be enough for the lower bound.

We consider N uniformly distributed random vectors $x_i\in \mathbb{S}^{n}$. Let $\varrho$ and $W$ be positive parameters, let $r=t_0+w$, where $w\in [-W,W]$.
For the purposes of the calculation we assume from the beginning that $W,\varrho\leq\frac{t_0}{20}$.
Consider a random polytope $Q$ in $\R^{n+1}$, defined as follows:
$$Q=\{x\in \R^{n+1} :\, \langle x,x_i\rangle\leq \varrho, \,\,\,\forall i=1,...,N\}.$$
Passing to the polar coordinates in $H_i=\{x:\,\langle x,x_i\rangle=\varrho\}$, we estimate the surface area of the half space $A_i=\{x:\,\langle x,x_i\rangle\leq\varrho\}$:
$$\gamma(\partial A_i)=\frac{1}{(n+1)\nu_{n+1}J_n}\int_{\R^{n}} e^{-\varphi(\sqrt{|y|^2+\varrho^2})}dy\succsim$$
$$\frac{1}{(n+1)\nu_{n+1} J_{n}} (n+1)\nu_n \int_{t_0-W}^{t_0+W} e^{-\varphi(t)} (t^2-\varrho^2)^{\frac{n-1}{2}}\frac{\sqrt{t^2-\varrho^2}}{t}dt\succsim$$
$$\frac{\sqrt{n}}{J_n}\left(1-\frac{\varrho^2}{(t_0-W)^2}\right)^{\frac{n}{2}} \int_{t_0-W}^{t_0+W} e^{-\varphi(t)} t^{n-1} dt.$$

Thus the expectation $\mathbb{E}\gamma(\partial Q)$ can be estimated from below by
\begin{equation}\label{expectation}
N\frac{\sqrt{n}}{J_n}\left(1-\frac{\varrho^2}{(t_0-W)^2}\right)^{\frac{n}{2}} \int_{t_0-W}^{t_0+W} e^{-\varphi(t)} t^{n-1} (1-p(t))^{N-1}dt,
\end{equation}
where $p(t)$ is the probability that the fixed point on the sphere of radius $t$ is separated from the origin by the hyperplane $H_i$.



As in \cite{fedia}, we use the formula for a surface area of a body of revolution to obtain the formula for $p(r)$:
\begin{equation}\label{uzhas}
p(r)=\left(\int_{-r}^{r}(1-\frac{t^2}{r^2})^{\frac{n-2}{2}}dt\right)^{-1}
\int_{\varrho}^{r}(1-\frac{t^2}{r^2})^{\frac{n-2}{2}}dt.
\end{equation}
By Laplace method, the first integral is approximately equal to $\frac{r}{\sqrt{n}}.$ Thus, after the change of variables $x=\frac{t}{r}$, we obtain
\begin{equation}\label{p(r)}
p(r)\approx \frac{\sqrt{n}}{r}r\int_{\frac{\varrho}{r}}^{1}(1-x^2)^{\frac{n-2}{2}}dx=\sqrt{n}\int_{\frac{\varrho}{r}}^{1}(1-x^2)^{\frac{n-2}{2}}dx.
\end{equation}
Notice, that for any $z\in(0,1)$,
\begin{equation}\label{el_ineq}
\int_{z}^{1} (1-t^2)^m dt\leq \frac{2}{z} \int_{0}^{1-z^2} s^m ds=\frac{2}{z(m+1)}(1-z^2)^{m+1}.
\end{equation}
By (\ref{p(r)}), (\ref{el_ineq}) applied with $z=\frac{\varrho}{r}$ and $m=\frac{n-2}{2}$, and the fact that $r\approx t_0,$

\begin{equation}\label{probability}
p(r)\precsim \frac{r}{\sqrt{n}\varrho} \left(1-\frac{\varrho^2}{r^2}\right)^{\frac{n}{2}}\precsim \frac{t_0}{\sqrt{n}\varrho} \left(1-\frac{\varrho^2}{(t_0+W)^2}\right)^{\frac{n}{2}}
\end{equation}
for all $r\in [t_0-W,t_0+W]$. At this point we choose
$$N=\frac{\sqrt{n}\varrho}{t_0} (1-\frac{\varrho^2}{(t_0+W)^2})^{-\frac{n}{2}}.$$
Observe that $(1-p(r))^{N-1}\lesssim (1-\frac{1}{N})^N\leq e^{-1}.$
Applying the above together with (\ref{expectation}) and (\ref{probability}), we get:
$$\mathbb{E}(\gamma(\partial Q))\succsim$$
\begin{equation}\label{expectation2}
\frac{n\varrho}{J_n t_0}\left(\frac{1-\frac{\varrho^2}{(t_0-W)^2}}{1-\frac{\varrho^2}{(t_0+W)^2}}\right)^{\frac{n}{2}} \int_{t_0-W}^{t_0+W} e^{-\varphi(t)} t^{n-1} dt.
\end{equation}


Let us now plug $W=\lambda t_0$. By Lemmas \ref{int_lam} and \ref{moments} and Remark \ref{to}, we observe, that
$J_{n-1}\approx \int_{t_0-\lambda t_0}^{t_0+\lambda t_0} t^{n-1} e^{-\varphi(t)}$. Thus,

$$\mathbb{E}(\gamma(\partial Q))\succsim\frac{n\varrho}{t_0}\times \frac{J_{n-1}}{J_n} \left(\frac{1-\frac{\varrho^2}{(t_0-W)^2}}{1-\frac{\varrho^2}{(t_0+W)^2}}\right)^{\frac{n}{2}} \approx  \frac{n\varrho}{t_0^2}\left(\frac{1-\frac{\varrho^2}{(t_0-W)^2}}{1-\frac{\varrho^2}{(t_0+W)^2}}\right)^{\frac{n}{2}}.$$

We plug $\varrho=\frac{1}{5\sqrt{\lambda n}}t_0$. Then
$$\frac{1-\frac{\varrho^2}{(t_0-W)^2}}{1-\frac{\varrho^2}{(t_0+W)^2}}\geq 1-\frac{1}{n},$$
which implies that
$$\mathbb{E}(\gamma(\partial Q))\succsim\frac{\sqrt{n}}{\sqrt{\lambda}t_0}.$$
This finishes the lower bound part of the Theorem \ref{main}.
$\square$


\section{Final remarks}
As was discussed in Section 3, Theorem \ref{main} entails Theorem \ref{th:main}. Its conclusion can be understood for any measure which has at least two bounded moments, so it is interesting to explore sufficiency of our conditions, i.e. spherical invariance and log concavity. We shall consider some examples of non rotation invariant or non log concave measures, for which the conclusion of Theorem \ref{th:main} does not hold.

\begin{example}\label{example1}
Consider Lebesgue measure concentrated on the cube $[-\frac{1}{2},\frac{1}{2}]^n$.
Due to convexity, the set of maximal surface area for this measure is the cube $[-\frac{1}{2},\frac{1}{2}]^n$ itself.
Its surface area is $2n$. However, $\mathbb{E}|X|\approx \sqrt{n}$ and $Var|X|\approx 1$ (see \cite{gian} for the proof), so if
Theorem \ref{th:main} was true, it would give $n^{\frac{1}{4}}$ as a maximal surface area. Thus there is no hope for Theorem \ref{th:main} to be true for all log concave measures. The isotropicity assumption would not change anything due to the homogeneity of Theorem \ref{th:main}.
\end{example}


\begin{example}\label{example3}
 Pick $\varepsilon\ll \frac{1}{n}$. We consider a rotation invariant non log concave measure $\gamma_{\varepsilon}$. Let its density be
$$
f(y) = c_n\begin{cases} 0 &\mbox{if } |y|\in [0,1-\varepsilon]\cup[1,\infty) \\
1 & \mbox{if } |y| \in(1-\varepsilon, 1). \end{cases}
$$

The normalizing constant
$$c_n=\nu_{n+1}(1-(1-\varepsilon)^{n+1})\leq (n+1)\varepsilon\nu_{n+1}.$$
For a random variable $X$ with density $f$ we compute
$$\mathbb{E}|X|=\frac{1-\frac{n+1}{2}\varepsilon}{1-\frac{n}{2}\varepsilon}+o((n+1)\varepsilon)\approx 1$$
and
$$Var |X|\approx\frac{\varepsilon^2}{4}.$$
Thus if Theorem \ref{th:main} was true, the maximal surface area would be of order $\frac{\sqrt{n}}{\sqrt{\varepsilon}}$. However,
$$\gamma_{\varepsilon}(\mathbb{S}^n)\geq \frac{1}{n\varepsilon^2},$$
which is greater than $\frac{\sqrt{n}}{\sqrt{\varepsilon}}$ for $\varepsilon\ll \frac{1}{n}$.
\end{example}

Example \ref{example3} shows, that for any dimension $n$ there exist a rotation invariant measure for which the conclusion of Theorem \ref{th:main} fails,
but it is hard to find an example of a density function which would serve all sufficiently large dimensions at once. It suggests the following conjecture.

\begin{conjecture}
Fix any real-valued function $\varphi(t)$ on the positive semi-axes. Then there exists a positive constant $C_{\varphi}$, depending on the function $\varphi(t)$,
such that for all $n\geq C_{\varphi}$,
$$\max_{Q\in \mathcal{K}_n}\gamma(\partial Q)\approx\frac{\sqrt{n}}{\sqrt{\mathbb{E}|X|} \sqrt[4]{Var|X|}},$$
where $X$ is a random vector on $\R^n$ distributed with respect to the density $e^{-\varphi(|X|)}$.
\end{conjecture}

\appendix

\section*{Appendix}

In this Appendix  we provide  a technical lemma which is believed to be well known to the specialists.
See \cite{kane} for the proof of the same statement in the case of Standard Gaussian Measure and polynomial level sets.

\begin{lemma}\label{convergence}
Let $\gamma$ be a probability measure on $\R^{n+1}$ with a continuous density $f(y)$. Then, for any convex set $Q$ in $\R^{n+1}$,
$$\int_{\partial Q} f(y) d\sigma(y)=\lim_{\varepsilon\rightarrow 0} \frac{\gamma(Q+\varepsilon B_2^{n+1})-\gamma(Q)}{\varepsilon},$$
where, as before, $d\sigma(y)$ stands for Lebesgue surface measure. 
\end{lemma}
\noindent\textbf{Proof.} For a convex set $Q$ in $\R^{n+1}$ and $\varepsilon>0$, we introduce the notation $A_{Q,\varepsilon}=\left(Q+\varepsilon B_2^{n+1}\right)\setminus Q$. We remark, that the normal vector $n_y$ is well defined almost everywhere for $y\in \partial Q$ if $Q$ is convex. So the function $f(y+tn_y)$ is defined almost everywhere on $\partial Q$. We shall apply the second Nazarov's system (\ref{second}), which we used in the proof of the main result. By convexity of $Q$,
\begin{equation}\label{outer_ann}
\gamma(A_{Q,\varepsilon})\geq \int_{\partial Q}\int_0^{\varepsilon} f(y+t\n_y)dtd\sigma(y),
\end{equation}
where the integration is understood in the Lebesgue sense. By Lebesgue Differentiation Theorem, for every $y\in\partial Q$ such that $n_y$ is defined,
$$\lim_{\varepsilon\rightarrow 0}\frac{1}{\varepsilon}\int_0^{\varepsilon} f(y+t\n_y)dt=f(y).$$
Consequently,
\begin{equation}\label{outer_ann_lim}
\lim_{\varepsilon\rightarrow 0}\frac{1}{\varepsilon}\gamma(A_{Q,\varepsilon})\geq \int_{\partial Q} f(y) d\sigma(y).
\end{equation}
On the other hand, we compare the measure of our annulus to the surface area of $Q+\varepsilon B_2^{n+1}$.

We note that for any $\varepsilon>0$ and $x\in A_{Q,\varepsilon}$ we may find $y\in \partial\left(Q+\varepsilon B_2^{n+1}\right)$ and $t\in [0,\varepsilon]$ so that $x=y-t\n_y$.

To see this, inscribe a ball centred at $x$ into $Q+\varepsilon B_2^{n+1}$ and chose $y$ to be a contact point of the ball and $\partial(Q+\varepsilon B_2^{n+1})$. We see, that $|x-y|\leq \varepsilon$, since
$$dist(x,\partial (Q+\varepsilon B_2^{n+1}))\leq \varepsilon.$$

We write
\begin{equation}\label{inner_ann}
\gamma(A_{Q,\varepsilon})\leq \int_{\partial (Q+\varepsilon B_2^{n+1})}\int_0^{\varepsilon} f(y-t\n_y)dtd\sigma(y).
\end{equation}
We observe, that
\begin{equation}\label{inner_ann_lim}
\lim_{\varepsilon\rightarrow 0}\frac{1}{\varepsilon}\gamma(A_{Q,\varepsilon})\leq \int_{\partial Q} f(y) d\sigma(y).
\end{equation}
Finally, (\ref{outer_ann_lim}) and (\ref{inner_ann_lim}) entail the conclusion of the Lemma. $\square$

\end{document}